%% file: manuscript.tex
\documentclass{llncs}

\usepackage{makeidx}
\makeindex

\usepackage[english]{babel}
\usepackage{tikz}
\usepackage{graphicx}
\usepackage{calc}
\usepackage{color}
\usetikzlibrary{automata,positioning,shapes,shadows,arrows}
\usepackage{amsfonts}
\usepackage{amsmath}

\usepackage{algpseudocode}
\usepackage{algorithm}
\usepackage{algorithmicx}
\usepackage{subcaption}
\title{Minimization and Equivalence in Multi-valued Logical Models of Regulatory Networks}
\author{Adam Streck, Therese Lorenz and Heike Siebert}
\institute{Freie Universit\"{a}t Berlin,
\email{adam.streck@fu-berlin.de}}
\date{\today}

\begin{document}

\maketitle

\begin{abstract}
Multi-valued logical models can be used to describe biological networks on a high level of abstraction based on the network structure and logical parameters capturing regulatory effects. Interestingly, the dynamics of two distinct models need not necessarily be different, which might hint at either only non-functional characteristics distinguishing the models or at different possible implementations for the same behaviour. Here, we study the conditions allowing for such effects by analysing classes of dynamically equivalent models and both structurally maximal and minimal representatives of such classes. Finally, we present an efficient algorithm that constructs a minimal representative of the respective class of a given multi-valued model.  

\end{abstract}
\begin{keywords}
regulatory networks, multi-valued models, minimization, equivalence, transition system
\end{keywords}
\section{Introduction}
Logical modelling has been used for describing and analysing biological systems such as gene regulatory networks for many years now, starting most notably with S.~Kauffman's work on Boolean networks~\cite{kauffman69a}. 
The resulting models can be viewed as non-uniform cellular automaton~\cite{dennunzio12a} with varying diameter. 
The underlying idea of capturing a system in terms of network structure and logical rules that govern the change of components has later been extended by R.~Thomas~\cite{thomas91a}. 
In his multi-valued logical framework components are allowed to adopt more than two activity levels and so-called logical parameters determine the effect of sets of regulators on a target component. 
The update scheme yielding transitions in the state space of the model can be chosen to allow for either all indicated value changes to be changed in one step, \emph{i.e.} synchronously, or for only one possible component update per step, \emph{i.e.} asynchronously. 
In the asynchronous update, as introduced by R. Thomas, this results in a non-deterministic transition system that generally is harder to analyse than in the synchronous, deterministic, case, but usually includes more realistic trajectories~\cite{thomas13a}.

In application, constructing a model, \emph{i.e.}, specifying the network structure and assigning the logical parameters, is a difficult task.
Usually not enough information is available to make unambiguous choices. Consequently, comparison of transition systems of many models becomes necessary. In this context, it is an interesting observation that different models may give rise to the same state transition system. 
For Boolean networks this phenomenon has been investigated and is directly related to superfluous edges--- edges without detectable dynamical effect---in the graph capturing the network structure (see \emph{e.g.} \cite{siebert08a}). 
In the multi-valued setting however, not only existence but also strength of a regulatory effect can be captured in the model. 
This allows for models differing both in structure and parameters to generate the same transition system, even if only functional edges are considered. 
This represent the full dynamical equivalence, as opposed to what is usually understood as dynamical equivalence in reduction studies, e.g.~\cite{naldi11a}, where the focus is only on the stabilizing behaviour only.

For a full understanding of this problem in the multi-valued case, the notion of functionality or \emph{observability} needs to be extended from the straight-forward intuition as it is usually used (see \emph{e.g.} \cite{streck12b}~\cite{gallet14a}). In this paper, we firstly present the proofs of results given originally in \cite{streck14a} that allow to describe the equivalence class of models generating the same transition system using a structurally maximal representative. This already allows to pinpoint the key properties to be focused on, namely self-loops and strength of effects as modelled by the value of a given logical parameter. However, the approach is not suited to application due to the high complexity. Extending the results in \cite{streck14a}, we then show in the second part of this paper that we can also systematically derive a structurally minimal representative of the equivalence class of a given model exploiting a resolved notion of observability and present an algorithm implementing our ideas.

\section{Background}
\input{example_toy_network}
We start by introducing the relevant notions, with a simple illustrative example given in Fig.~\ref{fig:toy_original}.

A multi-valued regulatory graph (RG) is a triple $G = (V,E,\rho)$ where:
\begin{itemize}
\item $V$ is a set of named \emph{components},
\item $E \subseteq V \times \mathbb{N}^+ \times V$ is a set of \emph{regulations},
\item $\rho : V \rightarrow \mathbb{N}^+$ assigns the maximal \emph{activity level} to a component.
\end{itemize}
It is required that for all $(u,n,v) \in E$ the non-zero integer $n$, called \textit{threshold}, is in the range of the source component \emph{i.e.} $0 < n \leq \rho(u)$.
We also use an additional function  $\theta : V \times V \to 2^{\mathbb{N}^+}$ which provides the \emph{thresholds} of all edges between two vertices and is defined as $\theta(u,v) = \{ n \mid (u,n,v) \in E \}$ where $u,v \in V$. 
Note that $\theta(u,v)=\emptyset$ if there are no edges from $u$ to $v$.

We use the symbol $\mathcal{G}$ to denote the set of multi-valued RGs. 
Also, in the following we use the notation $x_{i\gets n}$ for a substitution such that in a tuple $x$ the value indexed by $i$ is substituted for $n$.
\subsection{Discrete kinetic parameters}
\label{sec:discrete_kin_param}
We denote the \emph{state space} of a system with $S = \prod_{v \in V} [0,\rho(v)]$. Note that the state space is shared among the graphs that have the same function $\rho$ and thus also the same $V$. We will use  $\mathcal{G}_{\rho} = \{ (V,E,\rho') \mid \rho = \rho' \}$ to refer to the class of graphs that share the state space, \emph{e.g.}, the graphs in Figures~\ref{fig:toy_original},~\ref{fig:toy_reduced} belong to the same class.

The set $S$ represents all the qualitatively different configurations of a system. However, each component is dependent only on the values of its regulators. An equivalence class on $S$ w.r.t. regulation of a component $v \in V$ is called the regulatory context.
To define the relevant notions we first describe the \emph{activity interval} of a regulator. For formal reasons we consider an extended threshold function ${\Theta}$ with ${\Theta}(u,v) = \theta(u,v) \cup \{0, \rho(u) + 1\}$ for all $u,v \in V$. Then 
\begin{align}
I_v^u = \{[j,k) \mid  j,k \in {\Theta}(u,v) , j < k, \neg (\exists l \in {\Theta}(u,v) (j < l < k)) \}
\label{def:act_interval}
\end{align}
is the set of activity intervals of $u$ in regulation of $v$. Here, the intuition is that the regulator effect of $u$ on $v$ is constant in each of the intervals of $I^u_v$.
Note that $\bigcup I_v^u = [0, \rho(u)]$, even in the case that there are no edges from $u$ to $v$.

The set of regulatory contexts of $v$ is then denoted and defined $\Omega_v = \prod_{u \in V} I_v^u$. For each $v \in V$, a regulatory context $\omega \in \Omega_v$ is a product-of-sets where $\omega_u$ is the $u$-th set of the product for each $u \in V$. For convenience we use comparisons of activity levels and intervals. We say that $n < \omega_v, n \in [0,\rho(v)]$ \textit{iff} for each $m \in \omega_v$ we have $n < m$. Also note that the intervals are not overlapping, therefore if we extend $<$ to pair of intervals then $<$ naturally forms a total ordering.

The dynamics of the RG are given via integer values, called \emph{logical parameters}, assigned to each context.
We use a \emph{parametrization function} $K_v : \Omega_v \to [0, \rho(v)]$ for each $v \in V$.
The parametrization of a regulatory network $(V,E,\rho) \in \mathcal{G}$ is then a tuple $K = (K_1,\dots,K_{|V|})$. 
Lastly, we denote $\mathcal{K}_G$ the set of all parametrizations of the RG $G$, called the \emph{parametrization space}.

To denote whether an edge is superfluous in the system, we use the \textit{observability} constraint. 
Intuitively we say that an edge is not observable if the parameter value is the same for all pairs of regulatory contexts that differ only in presence of the said edge.
In Figure~\ref{fig:toy_original} it is apparent that the change in the $\omega_u$ does not cause a change in the parameter $K_v(\omega)$, \emph{i.e.}, the edge $(u,1,v)$ is not observable.

We now formalize this notion.
In the rest of the article we will use for an edge $(u,n,v)$ the symbol $n_- \in {\Theta}(u,v)$ for a value such that there is no $m \in {\Theta}(u,v)$ where $n_- < m < n$, and $n_+ \in {\Theta}(u,v)$ for a value such that there is no $m \in {\Theta}(u,v)$ where $n_+ > m > n$. 
Note that $n_-,n_+$ always exist as $0, (\rho(u) + 1) \in {\Theta}(u,v) - \theta(u,v)$.
Now, we say that $(u,n,v)$ is not observable in parametrization $K \in \mathcal{K}_G$ \textit{iff} for each $\omega \in \Omega_v$ such that $\omega_u = [n,n_+)$ we have $K_v(\omega) = K_v(\omega_{u \gets [n_-,n)})$.
  
\subsection{Asynchronous transition system}
Having a RG $G= (V,E,\rho)$ and a parametrization $K \in \mathcal{K}_G$ we can fully describe its dynamical behaviour as a \emph{transition system} (TS) over its state space $S$. 
This is a directed graph $(S,\rightarrow)$ where $\rightarrow \; \subseteq S \times S$ is the \emph{transition relation}. 
As mentioned, we are interested in asynchronous dynamics which means that the transition relation is non-deterministic.

First $K$ is converted into a so-called \emph{update function} $F^K=(F^K_v)_{v\in V}$ where $F^K_v : S \rightarrow [0,\rho(v)]$ for all $v \in V$. Here we exploit the fact that for each $s\in S$ and for each $v \in V$ there exists a context $\omega \in \Omega_v$ such that $s \in \prod_{u \in V} \omega_u$. To simplify the notation we will further write $s \in \omega$ instead of $s \in  \prod_{u \in V} \omega_u$. For every $v \in V$ we obtain the function $F^K : S \rightarrow S$ from a parametrization $K_v$ as
\begin{align}
F^K_v(s) = 
\begin{cases}
	s_v + 1, & \text{if } s_v < K_v(\omega), s \in \omega \text{,} \\
	s_v, & \text{if } s_v = K_v(\omega), s \in \omega \text{,} \\
	s_v - 1, & \text{if } s_v > K_v(\omega), s \in \omega \text{.}
\end{cases}
\label{def:update_fun}
\end{align}
Note that we have only three options how the value of a component can change, namely either increase by one, remain constant, or decrease by one. This provides us with a certain notion of direction of a derivative (positive, zero, or negative). In the following we will be using the term $\delta(F_v)(s) = F_v(s) - s_v$ to describe the partial derivative in the dimension $v$ in $s \in S$.

Having $F^K$, we now assign each parametrized RG a TS via the function $T_G : \mathcal{K}_G \to \{ (S,\rightarrow) \}$ where $T_G(K) = (S,\rightarrow)$ such that:
\begin{align} 
\forall s, s' \in S : s \rightarrow s' \iff (s \neq s') \wedge (F^K_v(s) = n) \wedge (\exists v \in V : s' = s_{v \gets n})\text{.}
\label{def:trans}
\end{align} 
In the following, we will compare the resulting TSs generated by RGs with the same state space, \emph{i.e.}, those in $\mathcal{G}_\rho$ for some $\rho$. 
 We denote this set of TSs $\mathcal{T}_{\rho} = \{ (S,\rightarrow) \mid T_G(K) = (S,\rightarrow), G \in \mathcal{G}_{\rho}, K \in \mathcal{K}_G\}$. 
An example of a shared TS is given in Figure~\ref{fig:toy_network}. 

\section{Conservative graph manipulations}
\label{sec:conservative_manipulations}
\input{example_ambiguity}
We now investigate the cases where different parametrizations of a single graph generate the same TS.
From (\ref{def:trans}) it is clear that two functions $F^{K} \neq F^{K'}$ will lead to distinct TSs, while coinciding functions $F^K, F^{K'}$ lead to the same dynamics.
We therefore focus on describing the situations where for $K \neq K'$ we still have $F^{K} = F^{K'}$.

Consider the simple example in Fig.~\ref{fig:ambiguity}. We see that putting $K'(\{0,1\},\{1,2\}) = 1$  instead $K(\{0,1\},\{1,2\}) = 2$ still yields $F^{K'}_v(0,1) = F^{K}_v(0,1) = 1$ and therefore the TS remains the same. 
This illustrates that, other than in the Boolean case, information on parameter values may get lost when deriving the update function. In Fig.~\ref{fig:ambiguous_source} we have a case where the parameter value lies outside its context and by incremental change we leave the context even before the value can be attained. 

We now define the notion of a \emph{canonical} parametrization that prohibits such effects. Observe that a value change in $v$ can cause the change of context only if $v$ regulates itself. Therefore we say that $K \in \mathcal{K}_G$ is \emph{canonical} if and only if
\begin{align}
\forall v \in V, \forall \omega \in \Omega_{v}, \omega_{v} = [j,k) : (K_v(\omega) \geq j - 1) \wedge (K_v(\omega) \leq k)\text{.}
\label{def:canonical}
\end{align}
We also denote $\mathcal{C}_G \subseteq \mathcal{K}_G$ the subset of canonical parametrizations in $\mathcal{K}_G$.

We can obtain a clear correspondence between $K$ and $F^K$ if all the contexts contain just a single state, so that no ambiguities are introduced in (\ref{def:update_fun}). This partition is achieved when only considering \emph{complete} graphs. For clarity we add that $(V,E,\rho) \in \mathcal{G}$ is complete if and only if for all $u,v \in V$ and every $n \in [1,\rho(u)]$ the edge $(u,n,v)$ is in $E$. This gives us the following theorem:
\begin{theorem}
\label{th:bijection}
For each $G = (V,E,\rho) \in \mathcal{G}$ it holds that if $G$ is complete then
$T_G$ defines a bijection between $\mathcal{C}_G$ and $\mathcal{T}_{\rho}$.
\end{theorem}
\begin{proof}
Let $T_G(K) = (S,\rightarrow)$ for some complete $G = (V,E,\rho)\in \mathcal{G}$ such that $K \in \mathcal{C}_G$. 
The one-to-one correspondence between $F^K$ and $(S,\rightarrow)$ immediately follows from (\ref{def:trans}). We therefore need to show that there is also such a correspondence between $K$ and $F^K$.

First, it is important to note that if the graph is complete, each regulatory context depicts only a single configuration. This is because if $G$ is complete then by (\ref{def:act_interval}) we have:
$$\forall v \in V : \Omega_v = \prod_{u \in V} \{[0,1),[1,2),\dots,[\rho(u),\rho(u)+1)\}\text{.}$$

Since each context is a singleton, each component has its value fixed. In that case, canonicity requires that the parameter in the context differs from the value only by $1$. 
Therefore we have only three options for the parameter value. 
More precisely, by substituting (\ref{def:canonical}) we have $$\forall v \in V, \forall s \in S\, :\, s_v-1 \leq K_v(\{s\}) \leq s_v + 1\text{.}$$

Then in such a case, (\ref{def:update_fun}) can be written as
\begin{align*}
F^K_v(s) = 
\begin{cases}
	s_v + 1, & \text{if } s_v + 1 = K_v(\{s\}) \text{,} \\
	s_v, & \text{if } s_v = K_v(\{s\})  \text{,} \\
	s_v - 1, & \text{if } s_v - 1 = K_v(\{s\})\text{,}
\end{cases}
\end{align*}
from which we immediately get $F^K_v(s) = K_v(\{s\})$. Thus (\ref{def:trans}) corresponds to:
\begin{align*} 
\forall v \in V, \forall s \in S : (K_v(\{s\}) = n) \wedge (K_v(\{s\}) \neq s_v) \iff s \rightarrow s_{v \gets n} \text{.}
\end{align*} 
\qed
\end{proof}
Based on this theorem, we can consider a complete graph with canonical parametrization as a representative of a class of models with the same behaviour. 
Now we show that it is possible to convert any graph with some parametrization into a complete graph with canonical parametrization, while keeping the dynamics unchanged.

First, we focus on the \emph{canonization} function $\textsf{Can} : \{ (G,K) \mid G \in \mathcal{G}, K \in \mathcal{K}_G \} \rightarrow \{ (G,C) \mid G \in \mathcal{G}, C \in \mathcal{C}_G \} $. For each component $v \in V$ and for each regulatory context $\omega \in \Omega_v$ with $\omega_v = [j,k)$ we construct $C$ as follows:
\[
C_{v}(\omega) = 
\begin{cases}
	j - 1, & \text{if } K_v(\omega) < j - 1 \\
	k, & \text{if } K_v(\omega) > k \qquad. \\
	K_{v}(\omega), & \text{otherwise } \\
\end{cases}
\]

The goal is to avoid that the parameter value cannot be reached in one transition from any state in the context.
The procedure is illustrated in the conversion from the parametrization $K$ in Figure~\ref{fig:ambiguous_source} to the parametrization $K'$ in Figure~\ref{fig:ambiguous_canonized}.  We now prove that this procedure indeed yields a canonical parametrization for any RG that shares the TS with the original one.
\begin{lemma}\label{th:canonization_cor}
Canonization is correct. For all $G \in \mathcal{G}$ and all $K \in \mathcal{K}_G$ it holds that if $\textsf{Can}(G,K) = (G,C)$ then $C$ is canonical.
\end{lemma}
\begin{proof}
There are two options for $K$ not to be canonical. The first option is that
$$\exists v \in V, \exists \omega \in \Omega_{v}, \omega_{v} = [j,k) : K_v(\omega) < j - 1$$

but then $C_v(\omega) = j - 1$, so $C$ is canonical. The second case  
$$\exists v \in V, \exists \omega \in \Omega_{v}, \omega_{v} = [j,k) : K_v(\omega) > k$$
can be treated analogously.
\qed 
\end{proof}
\begin{lemma}\label{th:canonization_cons}
Canonization is conservative. For all $G \in \mathcal{G}$ and all $K \in \mathcal{K}_G$ it holds that if $\textsf{Can}(G,K) = (G,C)$ then $T_G(K) = T_G(C)$.
\end{lemma}
\begin{proof}
Recall that the TSs $T_G(K)$ and $T_G(C)$ are fully defined by $F^K$ and $F^{C}$, respectively. We therefore need to show that $F^K = F^{C}$. 

For all $v \in V$ and for all $\omega \in \Omega_{v}$ the value $C_v(\omega)$ is set based on one of the three following cases. 

First consider the case that $K_v(\omega) < \omega_v$. 
Denote $\omega_v = [j,k)$. 
For all $s \in \omega$ it holds that $s_v > j - 1$ and then $C_v(\omega) = j - 1$. 
This means that for all $s \in \omega$ both  $K_v(\omega)$ and $C_v(\omega)$ are smaller than $s_v$ and therefore for each $s \in \omega$ we have $F^K_v(s)  = s_v - 1 = F^{C}_v(s)$.

The case that $K_v(\omega_v) > \omega_v$ can be treated analogously. 

The third case is that we have $K_{v}(\omega_v) = C_{v}(\omega_v)$ and thus by definition $F^K_v(s) = F^{C}_v(s)$ for any $s \in \omega_v$ .
\qed
\end{proof}
\begin{corollary}
For any $K$ holds: If $j = K_v(\omega)<\omega_v$ (resp.\ $j =K_v(\omega)>\omega_v$) then replacing $j$ with $j'$, $0\leq j'<\omega_v$ (resp.\ $\omega_v<j'\leq\rho(v)$) is conservative.
\label{th:param_manipulation_conservative}
\end{corollary}
Second, we extend the structure of a graph using the \emph{completion} function $\textsf{Comp}:\{ (G,K) \mid G \in \mathcal{G}, K \in \mathcal{K}_G \} \rightarrow \{ (\hat{G}, \hat{K}) \mid \hat{G} \in \mathcal{G}, \hat{K} \in \mathcal{K}_{\hat{G}} \}$. If  $G$ is complete, we map $(G,K)$ to itself.
For an incomplete $G =(V,E,\rho)$ and some $K \in \mathcal{K}_G$ we consider the non-empty set of missing edges $\hat{E} = \{(u,n,v) \mid u,v \in V, n \in [1,\rho(u)], (u,n,v) \notin E \}$. Assume that the set of all possible edges has some ordering. We extend the graph $G$ to $\hat{G}$ such that $\hat{G} = (V,E \cup \{min(\hat{E})\},\rho)$. 
As an example you can see that the RGs in Figure~\ref{fig:toy_original},~\ref{fig:toy_reduced} differ from each other only by the presence of a non-observable edge. Likewise for the GRs in Figure~\ref{fig:ambiguous_source},~\ref{fig:ambiguous_reduced}.

To extend the parametrization $\hat{K}$ to the new structure we first observe that $\hat{G}$ gives rise to new contexts that were obtained by partitioning some context of $G$ into two. 
To preserve the dynamical behaviour we simply assign the parameter value of the original context to both resulting new contexts. 
Have $(u,n,v) = min(\hat{E})$.
For each $v \in V$ and for each $\hat{\omega} \in \hat{\Omega_v}$ we then create $\hat{K}$ as:
\[
\hat{K}_{v}(\hat{\omega}) = 
\begin{cases}
	K_{v}(\hat{\omega}) & \text{if } v \neq v \vee (\omega_{u} = [j,k)  \wedge (j \neq n_{-} \vee k \neq n_{+}))\text{,} \\
	K_{v}(\hat{\omega}_{u \gets [ n_{-} , n_{+})}) & \text{otherwise.} \\
\end{cases}
\]
We now prove that for an incomplete RG we can use the completion procedure to add a new edge while retaining the dynamics.
\begin{lemma}
Completion is sound. For all $G \in \mathcal{G}$ and all $K \in \mathcal{K}_G$ it holds that if $\textsf{Comp}(G,K) = (\hat{G},\hat{K})$ then $\hat{G} \in \mathcal{G}$ and $\hat{K} \in \mathcal{K}_{\hat{G}}$.
\label{th:completion_sound}
\end{lemma}
\begin{proof}
We have $u,v \in V$ and $n \in \rho(u)$, therefore by definition of $\mathcal{G}$ we have that if $(V,E,\rho) \in \mathcal{G}$ then $(V,E \cup \{ (u,n,v) \},\rho) \in \mathcal{G}$.

The variables $n_-,n_+$ exist since $0$ is always a possible choice for $n_-$ and $\rho(u) + 1$ for $n_+$. For any $m \in \theta(u,v)$ we know that $(0 < m < \rho(u) + 1)$. 

From (\ref{def:act_interval}) we know that the only change occurs in the interval $I_{v}^{u}$. 
We therefore only need to show that $\hat{K}_{v}$ is extended to the affected contexts. Since $n_-,n_+$ exist, for any $\omega \in \Omega_{v}$ we have that $\omega_{u} \in \{i_1,\dots,[n_-,n_+),\dots,i_k\}$. Therefore for any $\hat{\omega} \in \hat{\Omega}_{v}$ we also have $\hat{\omega}_{u \gets [n_-,n_+)} \in \Omega_{v}$. Therefore $\hat{K}_v$ is defined on the whole $\hat{\Omega}_v$ for each $v \in V$.
\qed
\end{proof}
\begin{lemma}
Completion is conservative. For each $G \in \mathcal{G}$ and for each $K \in \mathcal{K}_G$ it holds that if $\textsf{Comp}(G,K) = (\hat{G},\hat{K})$ then $T_{G}(K) = T_{\hat{G}}(\hat{K})$.
\label{th:completion_conservative}
\end{lemma}
\begin{proof}
We have that $K$ differs from $\hat{K}$ only in a context $\hat{\omega} \in \Omega_{v}$ with $\hat{\omega} = [j,k)$ where either $j = n_-$ or $k = n_+$. 

Assume there is some $s \in \hat{\omega}$ for which $F^K_{v}(s) \neq F^{\hat{K}}_{v}(s)$. 
But we know that $[j,k) \subset [n_-,n_+)$ and therefore $s \in \hat{\omega}_{u \gets [ n_-, n_+)}$.
This would however imply that also $\hat{K}_{v}(\hat{\omega}_{v}) \neq  K_{v}(\hat{\omega}_{u \gets [ n_-, n_+)})$, which contradicts the definition of $\hat{K}$. 
We therefore have that $F^K(s) = F^{\hat{K}}(s)$.
\qed
\end{proof}
Since the completion procedure adds only one edge at a time, we need to repeat the procedure. This is captured in the following Lemma.
\begin{lemma}\label{th:completion_f}
For $G \in \mathcal{G}$, $K \in \mathcal{K}_G$, consider the
recursive sequence $\textsf{Comp}(G,K),$ $ \textsf{Comp}(\textsf{Comp}(G,K)), \dots$ \\
This sequence converges to a fixed point $(G^c,K^c)$  and $G^c$ is complete. 
\end{lemma}
\begin{proof}
The set $\hat{E}$ of missing edges in $G$ is finite as $V$ is finite. For each $v \in V$ also $[1,\rho(v)]$ is finite. In each iterative use of $\textsf{Comp}$ the size of $\hat{E}$ is decremented by one. The recursive sequence becomes constant when $\hat{E}$ is empty, indicating a fixed point $(G^c,K^c)$ of $\textsf{Comp}$. By definition, $G^c$ is complete.
\qed
\end{proof}
Combining all the statements above, we arrive at our final theorem:
\begin{theorem}\label{th:equivalence_comp}
Let $G,G' \in \mathcal{G}$, $K \in \mathcal{K}_G$, $K' \in \mathcal{K}_{G'}$ and denote $\textsf{Comp}^*(G,K)$, resp. $\textsf{Comp}^*(G',K')$ the fixed points derived from iterating $\textsf{Comp}$ starting in $(G,K)$ resp.\ $(G',K')$. \\
Then $T_G(K) = T_G'(K')$ \textit{iff} $\textsf{Can}(\textsf{Comp}^*(G,K)) = \textsf{Can}(\textsf{Comp}^*(G',K'))$.
\end{theorem}
\begin{proof}
We now know that $\textsf{Can}(\textsf{Comp}^*(G,K))$ and $\textsf{Can}(\textsf{Comp}^*(G',K'))$ are canonical and complete. The equivalence follows from $T_G$ being a bijection, as proven in Theorem~\ref{th:bijection}.
\qed
\end{proof}
\section{Network minimization}
Canonization and completion provide mathematical insights into the problem of dynamical equivalence. However in application it is more useful to have a minimal, rather than maximal structure. 
Also, while the canonization provides a good intuition about what the actual behaviour for each context is, it can have side-effects like converting a non-observable edge to observable, as is illustrated in Figure~\ref{fig:ambiguity}. 
We therefore introduce another form of parametrization, named \textit{normalized} parametrization, which prevents such effects but is more involved.
Using a \textit{normalization} procedure we then obtain a parametrization which is amenable to minimization. This section is divided into five consecutive steps:
\begin{enumerate}
\item We introduce a notion of \emph{observability} in the TS which allows us to see whether an edge is observable based on the transitions in the TS.
\item We introduce a notion of a \emph{monotone target value} (MTV) of a component. This value keeps the observability properties of TS, but is shared for a whole context.
\item We show how to compute the MTV from a parametrization and consequently how to compute a \textit{normalized} parametrization.
\item We show that every edge that is not observable in the TS is not observable in the respective normalized parametrization.
\item We introduce the minimization function for RGs based on normalized parametrization and explain how to test equivalence via minimization.
\end{enumerate}
\subsection{Observability in transition systems}
We have already defined the notion of observability in the parametrization. However we are mostly interested in the observability since it has implications on the dynamics. We now show how the it can be evaluated in the TS.

Intuitively, for an edge to be not observable there must be only a single value towards which the component evolves, no matter whether the regulator is above the thresholds of the said edge or below it. 
Formally for an RG $G = (V,E, \rho)$ is an edge $(u,n,v) \in E$ not observable in $(S,\rightarrow) = T_G(K)$ \textit{iff}:
\begin{eqnarray}
&\forall s \in S, s_u \in [n_{-}, n_{+}), \exists k \in [0, \rho(v)], \forall j \in [n_{-}, n_{+}) : \nonumber \\ 
&\delta(F_v^K)(s_{u \gets j}) = \textsf{Sgn}(k - (s_{u \gets j})_v)\,,
\label{def:ts_observability}
\end{eqnarray}
where $\textsf{Sgn} : \mathbb{Z} \rightarrow \{+1,0,-1\}$ is the usual sign function.

Consider the example in Figure~\ref{fig:toy_network}. It is easy to see that the regulation $(u,1,v)$ does not have any effect, since the left and right half of the TS are identical. Take in particular the example of the state $(v,u)=(0,0)$. We choose $k = 2$. Then for $j = 1$ it holds that $\delta(F_v^K)(0,1) = +1 = \textsf{Sgn}(2 - 0)$ and for $j = 0$ it holds that $\delta(F_v^K)(0,0) = +1 = \textsf{Sgn}(2 - 0)$. 

From the definition of the function F (\ref{def:update_fun}) and its derivative we easily see that:
\begin{eqnarray}
\label{th:equiv_delta_signum_parameter}
\forall s \in S, s_v\in\omega \in \Omega_v : \delta(F_v^K)(s) = \textsf{Sgn}(K_v(\omega) - s_v)\,.
\end{eqnarray}
This illustrates that the TS non-observability is the kind that we are interested in, since it relies on the actual dynamics of the network as captured in $F^K$. It is a stronger notion than the corresponding parametrization based one, since observability in the TS implies observability in the parametrization. We show the contraposition of this statement in the following lemma:
\begin{lemma}
\label{th:unobservable_parametrization}
Have $(u,n,v) \in E$ not observable in $K \in \mathcal{K}_G$. Then $(u,n,v)$ is not observable in $T_G(K)$.
\end{lemma}
\begin{proof}
By definition of observability (Section~\ref{sec:discrete_kin_param}) we have that for each $\omega \in \Omega_v$ such that $\omega_v = [n, n_{+})$ it holds that $K_v(\omega) = K_v(\omega_{u \gets [n_{-},n)})$. We therefore can set $k = K_v(\omega) = K_v(\omega_{u \gets [n_{-},n)})$ and from (\ref{th:equiv_delta_signum_parameter}) we immediately see that (\ref{def:ts_observability}) is satisfied for the whole range $[n_{-}, n_{+})$.
\qed
\end{proof}

\subsection{Monotone target value}
To relate the dynamics captured in a TS with a parametrization value, we introduce the notion of monotone target value. 
Intuitively, for a state $s \in S$ an MTV is a value towards which the component value $s_v$ evolves if we traverse only in the dimension of $v$ until $s_v$ either stabilizes or an opposite effect takes place. This idea is strongly linked to the derivative of the update function. 

Consider the example in Figure~\ref{fig:toy_dynamics}. Under the influence of the edge $(u,1,v)$ we see the trace $(0,1) \rightarrow (1,1) \leftrightarrow (2,1)$. Here the MTV for component $v$ in the state $(v,u)=(0,0)$ is $2$ since $\delta(F^K_v)(0,1) = +1 = \delta(F^K_v(1,1)) \neq \delta(F^K_v(2,1)) = -1$, \emph{i.e.}, at the level $s_v = 2$ an opposing effect takes place.

For a $T_G(K) = (S,\rightarrow)$, $v\in V$, $\omega\in\Omega_v$ and any $s \in \omega$ we denote the MTV by $(F_v^K)^{mon}(s)$ defined as:
\begin{eqnarray}
(F_v^K)^{mon}(s)=\begin{cases}
	s_v, & \text{if } \delta(F_v^K)(s)=0 \text{,} \\
	min\lbrace j>s_v\vert\delta(F_v^K)(s_{v \gets j})\neq +1 \rbrace, & \text{if } \delta(F_v^K)(s)=+1  \text{,} \\
	max\lbrace j<s_v\vert\delta(F_v^K)(s_{v \gets j})\neq -1 \rbrace, & \text{if } \delta(F_v^K)(s)=-1  \text{.}
\end{cases}
\label{def:monotone_update}
\end{eqnarray}
Also note that:
\begin{eqnarray}
\forall s\in S,\forall v\in V : \delta(F_v^K)(s)=\textsf{Sgn}((F_v^K)^{mon}(s)-s_v) \,.
\label{th:connection_delta_monotone_update}
\end{eqnarray}
The MTV can now be related to the observability in the TS. In particular, we can rewrite (\ref{def:ts_observability}) as:
\begin{eqnarray}
&\forall s \in S, s_u \in [n_{-}, n_{+}), \forall j \in [n_{-}, n_{+}) :  \nonumber \\ 
&\delta(F_v^K)(s_{u \gets j}) = \textsf{Sgn}((F_v^K)^{mon}(s) - (s_{u \gets j})_v)
\label{def:ts_mononte_obs}
\end{eqnarray}
\begin{lemma}
The non-observability conditions (\ref{def:ts_observability}) and (\ref{def:ts_mononte_obs}) are equivalent.
\end{lemma}
\begin{proof}
Clearly, (\ref{def:ts_mononte_obs}) implies (\ref{def:ts_observability}) since we can set $k = (F_v^K)^{mon}(s)$. 
For the other direction we show that if (\ref{def:ts_mononte_obs}) does not hold, then (\ref{def:ts_observability}) cannot hold either. 

If (\ref{def:ts_mononte_obs}) does not hold, then there is a state $s \in S$ and some $j \in  [n_{-}, n_{+})$  such that $\delta(F_v^K)(s') \neq \textsf{Sgn}((F_v^K)^{mon}(s) - s'_v)$ where $s'=s_{u \gets j}$. 
We distinguish three cases based on the value of $\delta(F_v^K)(s')$ and then again three cases based on the difference between $s_v$ and $s'_v$:
\vspace{10pt}\\
\textbf{Case $\delta(F_v^K)(s') = 0$ and $\textsf{Sgn}((F_v^K)^{mon}(s) - s'_v) \neq 0$:}
\begin{itemize}
\item If $s_v = s'_v$, then for any $k$ certainly $\textsf{Sgn}(k - s'_v) = \textsf{Sgn}(k - s_v)$. Also since $\textsf{Sgn}((F_v^K)^{mon}(s) - s'_v) \neq 0$ we get $(F_v^K)^{mon}(s) \neq s'_v=s_v$,  and thus $\delta(F_v^{K}(s)) \neq 0$ according to (\ref{th:connection_delta_monotone_update}). Then if there is a $k$ s.t. $\delta(F_v^K)(s') = \textsf{Sgn}(k - s'_v)$ then also $\textsf{Sgn}(k - s'_v) \neq \delta(F_v^{K}(s))$ and therefore (\ref{def:ts_observability}) does not hold.
\item If $s_v > s'_v$, it follows that $u = v$. From (\ref{def:monotone_update}) we get $(F_v^K)^{mon}(s) \geq s'_v$ since $\delta(F_v^K)(s') = 0$ indicating an effect change in the only problematic case that $\delta(F_v^K)(s) = -1$. Since $\textsf{Sgn}((F_v^K)^{mon}(s) - s'_v) \neq 0$ we have strict inequality $((F_v^K)^{mon}(s)) > s'_v$. Subsequently by (\ref{def:monotone_update}) there exists $l$ such that $s_v \geq l > s'_v$ and $\delta(F_v^{K})(s_{v \gets l}) \geq 0$. We have that $s'_v \in [n_{-}, n_{+})$ and $\delta(F_v^K)(s') = 0$, so to fulfil (\ref{def:ts_observability}) the $k$ must be chosen as $k = s'_v$. But $\delta(F_v^K)(s_{v \gets l}) \geq 0$ and $l \in (s'_v, s_v] \subseteq [n_{-}, n_{+})$ and therefore $k$ needs to satisfy $k > (s_{v \gets l})_v = l$. Together we get $k = s'_v < l \leq k$ which is a contradiction.
\item The case $s_v < s'_v$ can be treated analogously to $s_v > s'_v$.
\end{itemize}
\textbf{Case $\delta(F_v^K)(s') = +1$ and $\textsf{Sgn}((F_v^K)^{mon}(s) - s'_v) \leq 0$:}
\begin{itemize}
\item If $s_v = s'_v$, then, since $\textsf{Sgn}((F_v^K)^{mon}(s) - s'_v)\leq 0$ and thus $(F_v^K)^{mon}(s) \leq s'_v = s_v$, we have $\delta(F_v^K)(s) \leq 0$. Then $k$ in (\ref{def:ts_observability}) needs to satisfy $k \leq s_v$. Also $\delta(F_v^K)(s') = +1 = \textsf{Sgn}(k - s'_v)$ and thus $k$ needs to satisfy $k > s'_v$ leading to a contradiction.
\item The case $s_v > s'_v$ is impossible, since $s_v > s'_v \geq (F_v^K)^{mon}(s)$, so $\delta(F_v^K)(s) = -1$ and therefore by (\ref{def:monotone_update}) also $\delta(F_v^K)(s') = -1$ which is a contradiction.
\item If $s_v < s'_v$, then by (\ref{def:monotone_update}) there is $l$ such that $s_v \leq l < s'_v$ and $\delta(F_v^{K})(s_{v \gets l}) \leq 0$. As $l \in [n_{-}, n_{+})$, a $k$ satisfying condition (\ref{def:ts_observability})  must be such that $\textsf{Sgn}(k - (s_{v \gets l})_v) \leq 0$ and thus $k \leq l$. But $\delta(F_v^K)(s') = +1$ and a suitable $k$ must also satisfy $s'_v < k$. Together we again have the contradiction $k \leq l < s'_v < k$.
\end{itemize}
\textbf{Case $\delta(F_v^K)(s') = -1$ and $\textsf{Sgn}((F_v^K)^{mon}(s) - s'_v) \geq 0$:} \\
This case can be treated like the previous one.
\qed
\end{proof}

\subsection{Normalization algorithm}
We can see that the MTVs $(F_v^K)^{mon}(s)$ allow for a straightforward test of observability in the TS. Obtaining $(F_v^K)^{mon}(s)$ is however quite tedious. The size of $S$ is exponential w.r.t. the set $V$ and we have to unfold the TS to find the monotone paths characterizing the MTVs.
In this section we show that all states of a context share their MTV and additionally that we can obtain the MTV from a context directly.

We introduce the \textit{normalization} function, described in Algorithm~\ref{alg:norm}, that computes for each component and for each regulatory context of that component the MTV shared between the states of the context.
\begin{algorithm}                    
\caption{Calculate $\textsf{Norm}(K, v, \omega)$ where $\Theta(v,v) = \{n_0,...,n_k\}$. }      
\label{alg:norm}    
\begin{algorithmic}[1]
\State $[n_i, n_{i+1}) = \omega_v$
\If {$K_v(\omega) \in \omega_v$} 
    \State $\textsf{Norm}(K, v, \omega) = K_v(\omega)$
\ElsIf{$K_v(\omega) < \omega_v$}
	\State $\omega' = \omega_{v \gets [n_{i-1},n_{i})}$
    \If {$K_v(\omega') \geq n_i - 1$}
        \State $\textsf{Norm}(K, v, \omega) = n_i - 1$
    \Else 
    	\State $\textsf{Norm}(K, v, \omega) = \textsf{Norm}(K, v, \omega')$
    \EndIf
\Else 
	\State $\omega' = \omega_{v \gets [n_{i + 1},n_{i+2})}$
    \If {$K_v(\omega') \leq n_{i + 1}$}
        \State $\textsf{Norm}(K, v, \omega) = n_{i + 1}$
    \Else 
    	\State $\textsf{Norm}(K, v, \omega) = \textsf{Norm}(K, v, \omega')$
    \EndIf
\EndIf
\end{algorithmic}
\end{algorithm} 

In the algorithm we traverse through the contexts, rather than through states of a system, when looking for a monotone trace. As an example consider the $K'_v$ in Figure~\ref{fig:ambiguous_canonized} and $\omega = ([0,1),[1,2))$. Then $Norm(K',v,\omega) = Norm(K',v,\omega_{v \gets [n_{i+1},n_{i+2})}) = K'_v(\omega_{v \gets [n_{i+1},n_{i+2})}) = 2$. Note that this coincides with the value in Figure~\ref{fig:ambiguous_source}, which actually has a normalized parametrization.

The correctness of the approach is quite intuitive since a regulatory context is a subspace of the state space with uniquely determined target value.
This means that either the behaviour is monotone or there is exactly one stable state which breaks monotonicity in both directions. 
Considering the definitions of $(F_v^K)^{mon}(s)$ and the derivative. 
Easy calculations for the three cases $K_v(\omega)<\omega_v$,  $K_v(\omega)>\omega_v$ and  $K_v(\omega)\in\omega_v$ for a context $\omega$ immediately prove the following lemma.
\begin{lemma}
For any $v \in V$ and any $\omega \in \Omega_v$ we have $(F_v^K)^{mon}(s) = (F_v^K)^{mon}(s')$ for all $s,s' \in \omega$.
\label{th:monotonous_in_context}
\end{lemma}
Due to this Lemma we can extend the notion of MTV to regulatory contexts $\omega$ so that $(F_v^K)^{mon}(\omega) = (F_v^K)^{mon}(s)$ for any $s \in \omega$. Having this extension, we now prove the correctness of Algorithm~\ref{alg:norm}.
\begin{theorem}
For $v \in V$, $\omega \in \Omega_v$ it holds that $(F_v^K)^{mon}(\omega) = \textsf{Norm}(K, v, \omega)$.
\label{th:algorithm_correct}
\end{theorem}
\begin{proof}
Linking back to the importance of self-regulation already seen in the canonization, we lead the proof by induction w.r.t. the distance in number of activity intervals of self-regulation of $v$ between the context and its MTV.
We want to show that $\textsf{Norm}(K, v, \omega)$ returns the correct value after at most as many recursive calls as is the distance between the context and its MTV.
The notion of distance needed for this is defined as a function $D^K_v : \Omega_v \rightarrow \mathbb{N}$ with 
$$D^K_v(\omega) = \textsf{Max}(|\{ A\in I^v_v \mid \omega_v < A \leq A_\omega \}|, |\{ A\in I^v_v \mid \omega_v > A \geq A_\omega \}|)\text{,}$$ 
where $A_\omega\in I^v_v$ is the activity interval where it holds that $(F_v^K)^{mon}(\omega)\in A_\omega$.

Now we prove the theorem by the means of induction. Note that the proof does not constitute an invariant of the algorithm, as it proceeds in the other direction than the algorithm itself. In particular, we start by showing that for all the contexts that have their MTV within them or on their boundaries, the algorithm ends immediately with the correct value. Then we proceed to show that if the normalized parameter of any context whose distance is $m$ is correct and known, then the normalized parameter of a context whose distance is $m+1$ can be correctly determined by calling Algorithm~\ref{alg:norm} once.
\\
\\
\textbf{Base of induction (distance 0 and 1):}
\\
If $D^K_v(\omega) = 0$ then we know that $(F_v^K)^{mon}(\omega) \in \omega_v$. This implies that there is a state $s \in \omega$ such that $K_v(\omega) = s_v$  and $\delta(F^K_v)(s) = 0$ according to (\ref{def:update_fun}) and (\ref{def:monotone_update}). Therefore $(F_v^K)^{mon}(\omega) = K_v(\omega)= \textsf{Norm}(K, v, \omega)$, as set on the lines 2, 3. The recursion depth is 0.

If $D^K_v(\omega) = 1$, denote $\omega' = \omega_{v \gets A}$ such that $(F_v^K)^{mon}(\omega) \in A$. It follows from $D^K_v(\omega) = 1$ that $D^K_v(\omega') \leq 1$. We distinguish the two options:
\begin{itemize}
\item $D^K_v(\omega') = 0$: Then the above argument repeats and there is $s \in \omega$ such that $K_v(\omega') = s_v = (F_v^K)^{mon}(\omega')$. By the definition of the MTV we get $(F_v^K)^{mon}(\omega')=(F_v^K)^{mon}(\omega)$. Based on the ordering of $\omega,\omega'$ we arrive either on the line 9 or 16 of the algorithm and state correctly that $(F_v^K)^{mon}(\omega) = \textsf{Norm}(K, v, \omega) = \textsf{Norm}(K, v, \omega') = s_v $. 
The recursion depth is 1.
\item $D^K_v(\omega') = 1$: Since $D^K_v(\omega) = 1$ it follows that $(F_v^K)^{mon}(\omega')\in \omega_v$. In this case the respective MTVs take the adjacent values of the boundary between $\omega_v$ and $\omega'_v$. In the case that $\omega > \omega'$ we arrive on line 7 in the algorithm and assign $\textsf{Norm}(K, v, \omega) = n - 1$. This is correct as in any state of $\omega$ we monotonously update towards $\omega'$ and as we enter $\omega'$ by crossing the boundary value $n$, we change the direction back towards $\omega$, breaking the monotonicity. Analogously the correct value is assigned for the case $\omega < \omega'$. The recursion depth is 0.
\end{itemize}
\textbf{Induction step (distance over 1):}
\\
The induction assumption is that in at most recursion depth $m\geq 1$ the value of any $\omega' \in \Omega_v$ such that $D^K_v(\omega') \leq m$ is correctly set and consider now $D^K_v(\omega) = m+1$. 

In case $\omega_v > (F_v^K)^{mon}(\omega)$, since $D^K_v(\omega) > 1$, we have an $\omega'$ such that $D^K_v(\omega') = m\geq 1$ and $\omega_v > \omega'_v > (F_v^K)^{mon}(\omega)$. 
From the definition of the MTV it follows that all states in both $\omega$ and $\omega'$ monotonously decrease under $F_v^K$, therefore also $ \omega'_v > (F_v^K)^{mon}(\omega')$ which gives us $(F_v^K)^{mon}(\omega) = (F_v^K)^{mon}(\omega')$. In the algorithm this is assured on line 9 and by induction hypothesis $(F_v^K)^{mon}(\omega')$ is correctly determined by the algorithm in at most $m$ recursions, giving us the desired result for $(F_v^K)^{mon}(\omega)$.
The case that $\omega_v < (F_v^K)^{mon}(\omega)$ is again analogous, leading to line 14 instead of 9. The recursion depth is now $m+1$.

Since the set of intervals is finite and the recursion traverses mo\-no\-to\-nous\-ly, we terminate in the recursion depth of at most $\textsf{Max}(\{ \rho(v) | v \in V \})$.
\qed
\end{proof}
Using the normalization function we can, similarly to canonization, create a conservative and sound transformer on parametrizations. We extend $\textsf{Norm}$ to a function $\textsf{Norm} : \{ (G,K) \mid G \in \mathcal{G}, K \in \mathcal{K}_G \} \rightarrow \{ (G, N) \mid G \in \mathcal{G}, N \in \mathcal{N}_G \} $ where $\mathcal{N}_G  \subseteq \mathcal{K}_G$ is the set of normalized parametrizations of $G = (V,E,\rho)$ and $\textsf{Norm}(G,K)= (G,N)$ where $N$ is defined by
$$\forall v \in V, \forall \omega \in \Omega_v : N_v(\omega) = \textsf{Norm}(K,v, \omega)\text{.}$$
We have proven correctness of normalization already in Theorem~\ref{th:algorithm_correct}, so it only remains to prove that normalization is conservative.
\begin{lemma}
Normalization is conservative. For all $G \in \mathcal{G}$ and every $K \in \mathcal{K}_G$ it holds that if $\textsf{Norm}(G,K) = (G,N)$ then $T_G(K) = T_G(N)$.
\label{th:nat_correct}
\end{lemma}
\begin{proof}
Observe that if $K_v(\omega) \in \omega_v$ then $\textsf{Norm}(K, v,\omega) = K_v(\omega)$. In Corollary~\ref{th:param_manipulation_conservative} we have shown that if $K_v(\omega) < \omega_v$, then it is conservative to replace $K_v(\omega)$ with any $l \in \mathbb{N}$ such that $l < \omega_v$, which is also the case in Algorithm~\ref{alg:norm}. The same holds for the case that $K_v(\omega) > \omega_v$. \qed
\end{proof}
\subsection{Observability in normalized parametrization}
We have seen now that observability in the sense of an actual dynamical effect should not be evaluated based on the parametrization but rather on the TS. All information needed to construct a TS is captured in the MTVs due to its relation to the derivative and thus the update function. At the same time, an important aspect of parametrizations is shared, namely that the MTV stays fixed within a context. 
This allows us to link observability in parametrization and TS, as is shown in the following theorem that complements Lemma 6.
\begin{theorem}
For all $G \in \mathcal{G}$ and every $K \in \mathcal{K}_G$ it holds that if $\textsf{Norm}(G,K) = (G,N)$ then every edge that is not observable in $T_G(K)$ is not observable in $N$.
\label{th:normalized_is_unobservable}
\end{theorem}
\begin{proof} Assume that the above does not hold, \emph{i.e.} there exists an edge $(u,n,v) \in E$ s.t. (\ref{def:ts_mononte_obs}) holds, but also it holds that:
\begin{eqnarray}
&\exists \omega \in \Omega_v, \omega_u = [n,n_{+}), \omega^\downarrow = \omega_{u \gets [n_{-},n)}  : N_v(\omega) \neq N_v(\omega^\downarrow)\text{.}
\label{def:contradiction_to_obs}
\end{eqnarray}
\textbf{Case $u \neq v$:}\\
Note that in this case we have $\omega_v = \omega^\downarrow_v$.

We have $\delta(F^{N}_v)(s)=\delta(F^{K}_v)(s)$ for all states $s \in \omega$ as can be easily deduced from Lemma~\ref{th:nat_correct}. For all $s\in\omega, s'\in\omega^\downarrow$ we have $(F_v^K)^{mon}(s)=N_v(\omega)\neq N_v(\omega^\downarrow)= (F_v^K)^{mon}(s')$. It follows that we can only meet the condition $\textsf{Sgn}((F_v^K)^{mon}(s) - s_v)=\textsf{Sgn}((F_v^K)^{mon}(s) - s'_v)=\delta(F^{K}_v)(s') $ for all $s\in\omega,s'=s_{u \gets j}, j\in\omega^\downarrow$ as demanded in (\ref{def:ts_mononte_obs}) if and only if either $N_v(\omega) < \omega_v \wedge N_v(\omega^\downarrow) < \omega^\downarrow_v$ or $N_v(\omega) > \omega_v \wedge N_v(\omega^\downarrow) > \omega^\downarrow_v$.

First consider that $(N_v(\omega) < \omega_v) \wedge (N_v(\omega^\downarrow) < \omega^\downarrow_v)$. 
If the condition on the line 7 is satisfied for both $\omega$ and $\omega^\downarrow$ we immediately see that $N_v(\omega) = N_v(\omega^\downarrow)$.
If it is satisfied for exactly one, then apparently we break the requirement (\ref{def:ts_mononte_obs}) as for $s \in \omega'_v$ we have $(s_{u \gets n_-})_v = s_v$ but $\delta(F^K_v)(s_{u \gets n_-}) \neq \delta(F^K_v)(s)$. 

We therefore meet the condition on the line 8 and from the line 9 we know that for $\omega'$ as defined there $N_v(\omega') = N_v(\omega) \neq N_v(\omega^\downarrow) = N_v((\omega^\downarrow)')$.
Since  $N_v(\omega') \neq N_v((\omega^\downarrow)')$ we have again that $N_v(\omega') \not \in \omega'_v$ and the same for $(\omega^\downarrow)'$.
Therefore it again must hold that $N_v(\omega') < \omega'_v$ and $N_v((\omega^\downarrow)') < (\omega^\downarrow)'_v$. 
Apparently, the argument is recursive, requiring that for each $\omega' \in \Omega_v$ such that $\omega'_v < \omega_v$ it holds that $N_v(\omega') < \omega'_v$.
But then ultimately $N_v(\omega) < 0$ which contradicts the definition of $K$.

For the case that $N_v(\omega) > \omega_v$ a similar argument holds using the upper boundary $N_v(\omega) \leq \rho(v)$.
\\
\textbf{Case $u = v$:}\\
First note that in this case we have $\omega' = \omega^\downarrow$ for $\omega'$ as defined in the algorithm. In case that $k =(F^K_v)^{mon}(s)<n$ for any $s \in \omega$ we execute $\textsf{Norm}(K,v,\omega)$ in the algorithm and the condition on the line 4 is met. Then depending on the value $K(\omega')$:
\begin{itemize}
\item $K(\omega') \geq n - 1$: then $k = n - 1$, meaning that for any $s' \in \omega'$ we have $(F^K_v)^{mon}(s')_v = n - 1$, otherwise the edge would be observable. But then also $N_v(\omega) = N_v(\omega')$ which contradicts (\ref{def:contradiction_to_obs}).
\item $K(\omega') < n - 1$: then according to the line 9 we have $\textsf{Norm}(K, v, \omega) = \textsf{Norm}(K, v, \omega')$, again contradicting (\ref{def:contradiction_to_obs}).
\end{itemize}
The case that $k \geq n$ can be treated similarly.
\qed 
\end{proof}

\subsection{Minimizing the model}
Since we know that after normalization, all non-observable edges can be directly detected, constructing a reduction algorithm is rather straight-forward. Our minimization process is closely related to the completion process in Section~\ref{sec:conservative_manipulations}.

We use the \emph{minimization} function $\textsf{Min}:\{ (G,N) \mid G \in \mathcal{G}, N \in \mathcal{N}_G \} \rightarrow \{ (\check{G},\check{N}) \mid \check{G} \in \mathcal{G}, \check{N} \in \mathcal{N}_{\check{G}} \}$ to eliminate the non-observable edges.
If $(G, N)$ is minimal, \emph{i.e.} there are no non-observable edges, we map it to itself. Otherwise have $\check{E}$ the set of non-observable edges in $G = (V,E,\rho)$ and an arbitrary total ordering on $\check{E}$ and $(u, n, v) = min(\check{E})$, then:
\begin{enumerate}
\item $\check{G} = (V, E \setminus \{ (u, n, v) \},\rho)$,
\item $\check{N}_{v}(\omega_{u \gets [n_{-},n_{+})}) = N_{v}(\omega_{u \gets [n_{-},n)}) = N_{v}(\omega_{u \gets [n,n_{+})})$.
\end{enumerate}
The nature of MTVs ensures that $\check{N}$ is again in $\mathcal{N}_{\check{G}}$.
As an example consider the network in Figure~\ref{fig:ambiguous_source}. The edge $(v,1,v)$ is apparently not observable. We therefore remove it from $E$ and with that we set $\check{N}_{v}([0,3),[0,1)) = N_{v}([0,1),[0,1))  = N_{v}([1,3),[0,1)) = 0$ and $\check{N}_{v}([0,3),[1,2)) = N_{v}([0,1),[1,2))  = N_{v}([1,3),[1,2)) = 2$. Note that this coincides with the network in Figure~\ref{fig:ambiguous_reduced}, which is in fact minimized and normalized.

This procedure can be seen as an inversion to the completion as defined in Section~\ref{sec:conservative_manipulations},
where we created two new contexts by splitting one, keeping the values, whereas here we merge two contexts with the same value into one.
Note that the fixed points of $\textsf{Comp}$ and $\textsf{Min}$ are not dependent on the order on  $V \times \mathbb{N} \times V$. We can therefore take arbitrary, but fixed, order and execute both $\textsf{Comp}$ and $\textsf{Min}$ according to this order. Then
for $\textsf{Min}(G, N) = (\check{G}, \check{N})$ with $N \in \mathcal{N}$ and $G \neq \check{G}$ we have $\textsf{Comp}(\check{G}, \check{N}) = (G, N)$. This can be easily verified by applying the two operations successively.
Since completion is sound and conservative, minimization in such corresponding cases is also sound and conservative. 
The only remaining case is that $\textsf{Min}(G, N) = (G, N)$, but this is obviously sound and conservative too.

Iteration of the $\textsf{Min}$ function will then lead to a minimal structure, as demonstrated by the following lemma:
\begin{lemma}
For $G \in \mathcal{G}$, $K \in \mathcal{K}_G$ and $\textsf{Norm}(G,K) = (G,N)$ consider the
recursive sequence $\textsf{Min}(G,N),$ $ \textsf{Min}(\textsf{Min}(G,N)), \dots$ \\
 This sequence converges to a fixed point $(G^m,N^m)$  and $G^m$ is minimal, \emph{i.e.} there are no other $(G',K')$ s.t. $T_G(K) = T_{G'}(K')$ and $G'$ has less edges than $G^m$.
\label{th:minimization_full}
\end{lemma}
\begin{proof}
Since we change $N$ only if we remove an edge and the edge set is finite, the existence of a fixed point is trivial. If $G^m$ was not minimal, then inevitably there would have to be an edge which is not observable in $T_G(K)$ but observable in $N$ which contradicts Theorem~\ref{th:normalized_is_unobservable}.
\qed
\end{proof}
Now we can conclude the section with a theorem about equivalence checking through minimization, complementing the result of Section~\ref{sec:conservative_manipulations} and showing that the set of dynamically equivalent RGs has both maximal and minimal elements w.r.t. set inclusion on the regulators.
\begin{theorem}\label{th:equivalence_min}
Let $G,G' \in \mathcal{G}$, $K \in \mathcal{K}_G$, $K' \in \mathcal{K}_{G'}$ and denote $\textsf{Min}^*(\textsf{Norm}(G,K))$ and $\textsf{Min}^*(\textsf{Norm}(G',K'))$ the fixed points derived from iterating $\textsf{Min}$ starting in $\textsf{Norm}(G,K)$ resp.\ $\textsf{Norm}(G',K')$. \\ Then $T_G(K) = T_G'(K')$ \textit{iff} $\textsf{Min}^*(\textsf{Norm}(G,K)) = \textsf{Min}^*(\textsf{Norm}(G',K'))$.
\end{theorem}
\begin{proof}
As in the proof of Theorem~\ref{th:equivalence_comp} we now know that both $\textsf{Min}^*(\textsf{Norm}(G,K))$ and $\textsf{Min}^*(\textsf{Norm}(G',K'))$ are normalized and minimal. Due to Theorem~\ref{th:normalized_is_unobservable} we know that in the minimization exactly the non-observable edges in the TS are removed.
So if the two TSs are equivalent, the set of remaining edges is the same in both $G$ and $G'$. 

Additionally we know that the parameter values are not changed during the minimization, only in the normalization where the value is set to the MTV of the states covered by the context.
The only way how we could have a difference between $\textsf{Min}^*(\textsf{Norm}(G,K))$ and $\textsf{Min}^*(\textsf{Norm}(G',K'))$ is that there is a state $s \in T_G(K)$ and a component $v$ such that $(F^K_v)^{mon}(s) \neq (F^{K'}_v)^{mon}(s)$. 
But then by (\ref{def:monotone_update}) we have some $r \in T_G(K)$ such that $\delta(F^K_v)(r) \neq \delta(F^{K'}_v)(r)$ and $T_G(K) \neq T_{G'}(K')$.
\qed
\end{proof}

\section{Complexity}
It has been shown that testing equivalence for two expressions over $n$ variables is co-NP complete~\cite{bloniarz84a}. Since in a Boolean network a component $v \in V$ can be regulated by up to $|V|$ nodes, the update function can be an expression over $|V|$ variables and deciding whether two parametrizations are equivalent is therefore necessarily exponential w.r.t. the number of components. However in the minimization approach we can parametrize by the in-degree. For a given  RG $G = (V,E,\rho)$, denote $v^- = \{ (u,n,v) | \exists n \in \mathbb{N}, \exists u \in V, (u,n,v) \in E \}$ and set $k = \textsf{Max} (\{ |v^-| \mid v \in V \})$, then:
\begin{theorem}
$\text{TIME}(\textsf{Min}^*(\textsf{Norm}(G,K))) = \mathcal{O}(|V| \cdot 2
^{k} \cdot k)$.
\label{th:complexity}
\end{theorem}
\begin{proof}
There are at most $2^k$ different regulatory contexts for any $v \in V$. As we call $\textsf{Norm}(K, v, \omega)$ for each $\omega \in \Omega_v$ no recursion is needed---either we set the value directly or we set it equal to some other value that will eventually be known. We therefore call $\textsf{Norm}(K,v,\omega)$ at most $2^{k}$ times for each $v \in V$.

In the minimization part we need to obtain the set $\check{E}$. 
For an edge $(u,n,v) \in E$
to detect whether $(u,n,v) \in \check{E}$ we need to consider all $ \omega \in \Omega_v$ with  $\omega_u = [n,n_{+})$ and compare them to the respective $\omega_{u \gets [n_{-}, n)}$. At most we need to do $(|\Omega_v| / 2) \in \mathcal{O}(2^k)$ pair-wise comparisons. Since we have to consider at most $k$ regulators of $|V|$ components, we obtain in total $\mathcal{O}(|V| \cdot 2
^{k} \cdot k)$. To remove the edge we remove the set of tested contexts, which again leads to $\mathcal{O}(|V| \cdot 2
^{k} \cdot k)$. 
\qed
\end{proof}
Intuitively, since we only change the values in place or remove them, no additional space is needed. The spatial complexity is thus equal to the size of the input.
\section{Changing the update scheme}
Throughout the article we have been considering only the asynchronous update function, which was advantageous for the definition of non-observability in the TS and the monotone update. 
However we can extend our result to arbitrary update schemes where $F_v(s)$ can be applied for any $v \in V$ and $s \in S$, see~\cite{gershenson04a} for examples.
As an example, we consider the synchronous update. 
Then $T_G(K) = (S,\rightarrow)$ such that:
$$\forall s, s' \in S : s \rightarrow s' \iff \forall v \in V : F_v(s) = s'_v\text{.}$$
We can therefore immediately see that two TSs are again identical \textit{iff} the respective update functions are identical. 
Since the conversion from parametrization to update function remains unchanged, we have that two parametrizations are dynamically equivalent under the asynchronous update scheme \textit{iff} they are dynamically equivalent under any other scheme that depends on $F$.

\section{Conclusion}
In the article we have thoroughly investigated the notion of observability for multi-valued logical models of regulatory networks with an arbitrary update function.
We have shown that, unlike in the case of Boolean networks, a mere test on equality of parameters in contexts with and without an edge is not sufficient to determine if the edge has a detectable impact on the network dynamics. As illustrated, at the heart of such an effect lies the connection of self-regulation and strength of regulation in general, as encoded in the parameter values. Here, self-regulations can act as amplifiers for the effect of other regulators. Consequently the same effect could be achieved without such a self-regulation simply by strengthening the original regulatory effect, \emph{i.e.} by adapting the corresponding parameter value.

One way to obtain an unambiguous representation of a given model that characterizes the whole class of dynamically equivalent models utilizes the notion of canonical parametrization. Here, in essence, the effects of regulations are standardized, but to do so, the underlying structure, \emph{i.e.}, the set of regulators for each component, has to be blown up. Although fully answering the mathematical problem of characterizing the equivalence classes, this approach suffers from high computational complexity in practice.

Exploiting a refined notion of edge observability directly tied to the dynamical behaviour as captured in the transition system, we were able to construct a minimal
parametrized RG w.r.t. the set inclusion on the set of regulators as representative of a dynamical equivalence class. The key step to make the construction feasible for application was finding a parametrization where the refined and the original notion of observability coincides. We presented an efficient and simple algorithm that allows us to obtain this normalized parametrization from an arbitrary one.

Our thorough mathematical investigation uncovered a minimal representative that is helpful not only in model comparison when modelling under uncertainties but also for understanding the regulatory effects underlying an observed behaviour. However, it should be mentioned that other models in the class might also carry interesting information. In essence they represent the different possible mechanisms that can be used to implement a desired effect. In particular, in biological applications comparison between the actual and the possible implementations might help to gain a deeper understanding about the constraints governing biological systems. 

From the practical perspective, since the algorithm has very low complexity, we would suggest employing it in tools that focus on parametrization space analysis to prevent costly computations where they are not needed.
For future work we propose to attempt to translate the solution from the multi-valued to piece-wise linear models, which share the threshold behaviour, but usually feature even larger parametrization space and therefore would probably benefit even more from reduction techniques.
\bibliographystyle{splncs}
\bibliography{manuscript}
\end{document}

%% file: example_toy_network.tex
\begin{figure}[t]
	\centering
	\renewcommand{\arraystretch}{1.33}
	\begin{subfigure}[b]{0.3\textwidth}
		\centering
		\begin{tikzpicture}
	      \node[draw, circle, thick, minimum size=8mm] (1) at (0,0) {$v$};
	      \node[draw, circle, thick, minimum size=8mm] (2) at (1.5,0) {$u$};
	      \draw [->, in=105, out=75, loop] (1) to node[above] {$2$} (1) ;
	      \draw [->, bend left=20] (2) to node[above] {$1$} (1) ;
	      \draw [->, bend left=20] (1) to node[above] {$1$} (2) ;
	    \end{tikzpicture}\\
    	$V = \{ v , u\} $\\
		$E = \{ (u,1,v),$  $(v,2,v), (v,1,u) \}$\\
		$\rho(v) = 2, \rho(u)=1$ \\
    	\begin{tabular}{|c|c|}
    	\hline
    	$(\omega_v, \omega_u)  \in \Omega_v$ & $K_v(\omega)$  \\ \hline
    	$([0,2),[0,1))$ & 2 \\
    	$([2,3),[0,1))$ & 1 \\
      	$([0,2),[1,2))$ & 2 \\
    	$([2,3),[1,2))$ & 1 \\  	
    	\hline
    	\end{tabular}
    	\begin{tabular}{|c|c|}
    	\hline
    	$(\omega_v, \omega_u)  \in \Omega_u$ & $K_u(\omega)$  \\ \hline
    	$([0,1),[0,2))$ & 0\\
    	$([1,3),[0,2))$ & 1\\
    	\hline
    	\end{tabular}
	    \caption{}
	    \label{fig:toy_original}
    \end{subfigure}
	\begin{subfigure}[b]{0.025\textwidth}
	$\rightarrow$\vspace{80pt}
    \end{subfigure}
    \begin{subfigure}[b]{0.3\textwidth}
    \centering
    \begin{tikzpicture}
      \node[draw, rectangle, thick, minimum width=10mm, minimum height=5mm] (00) at (0,0) {(0,0)};
      \node[draw, rectangle, thick, minimum width=10mm, minimum height=5mm] (10) at (0,-1) {(1,0)};
      \node[draw, rectangle, thick, minimum width=10mm, minimum height=5mm] (20) at (0,-2) {(2,0)};
      \node[draw, rectangle, thick, minimum width=10mm, minimum height=5mm] (01) at (1.25,0) {(0,1)};
      \node[draw, rectangle, thick, minimum width=10mm, minimum height=5mm] (11) at (1.25,-1) {(1,1)};
      \node[draw, rectangle, thick, minimum width=10mm, minimum height=5mm] (21) at (1.25,-2) {(2,1)};
      \draw [->, thick] (01) to (11);
      \draw [->, thick] (11) to (21);
      \draw [->, thick] (21) to (11);
      \draw [->, thick] (00) to (10);
      \draw [->, thick] (10) to (20);
      \draw [->, thick] (20) to (10);
      
      \draw [->, thick] (01) to (00);
      \draw [->, thick] (10) to (11);
      \draw [->, thick] (20) to (21);
    \end{tikzpicture}
    \vspace{3pt}
		\caption{}
	    \label{fig:toy_dynamics}
    \end{subfigure}
	\begin{subfigure}[b]{0.025\textwidth}
	$\leftarrow$\vspace{80pt}
    \end{subfigure}    
    \begin{subfigure}[b]{0.3\textwidth}
		\centering
		\begin{tikzpicture}
	      \node[draw, circle, thick, minimum size=8mm] (1) at (0,0) {$v$};
	      \node[draw, circle, thick, minimum size=8mm] (2) at (1.5,0) {$u$};
	      \draw [->, in=105, out=75, loop] (1) to node[above] {$2$} (1) ;
	      \draw [->] (1) to node[above] {$1$} (2) ;
	    \end{tikzpicture}\\
    	$V' = \{ v , u\} $\\
		$E' = \{ (v,2,v), (v,1,u)  \}$\\
		$\rho(v) = 2, \rho(u)=1$ \\
    	\begin{tabular}{|c|c|}
    	\hline
    	$(\omega_v, \omega_u) \in \Omega'_v$ & $K'_v(\omega)$ \\ \hline
   		$([0,2),[0,2))$ & 2 \\
    	$([2,3),[0,2))$ & 1 \\
    	\hline
    	\end{tabular}
    	\begin{tabular}{|c|c|}
    	\hline
    	$(\omega_v, \omega_u)  \in \Omega'_u$ & $K'_u(\omega)$  \\ \hline
    	$([0,1),[0,2))$ & 0\\
    	$([1,3),[0,2))$ & 1\\
    	\hline
    	\end{tabular}
	    \caption{}
	    \label{fig:toy_reduced}
    \end{subfigure}
	\caption{ Toy example. a) A network $G$ with the parametrization $K$. c) A reduced toy network $G'$ with the parametrization $K'$. Note that for any $x \in [0, \rho(v)]$ and $y \in [0, \rho(u)]$ it holds that $F^K_v(x,y) = F^{K'}_v(x,y)$ and $F^K_u(x,y) = F^{K'}_u(x,y)$. Therefore we have b) the transition system $T_G(K) = T_G'(K')$.}
	\label{fig:toy_network}
\end{figure}
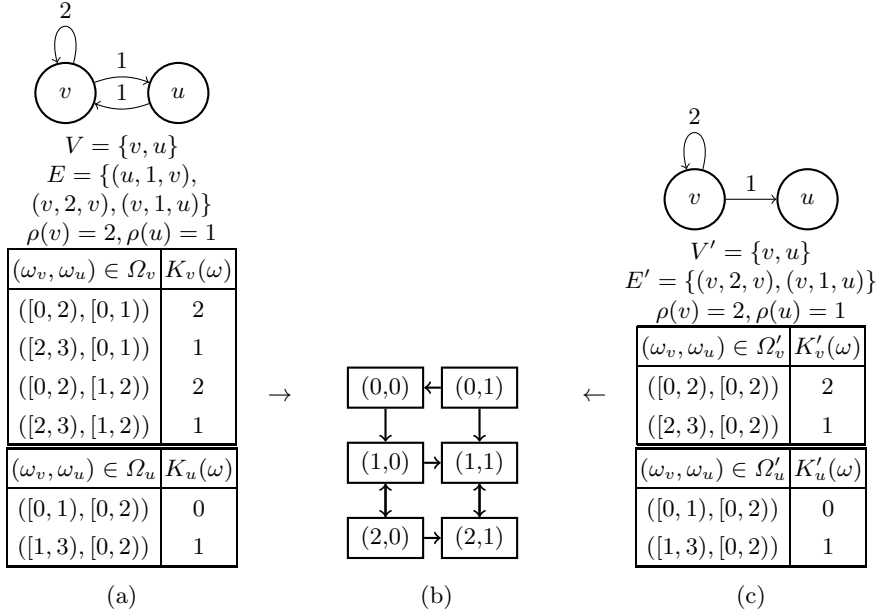

%% file: example_ambiguity.tex
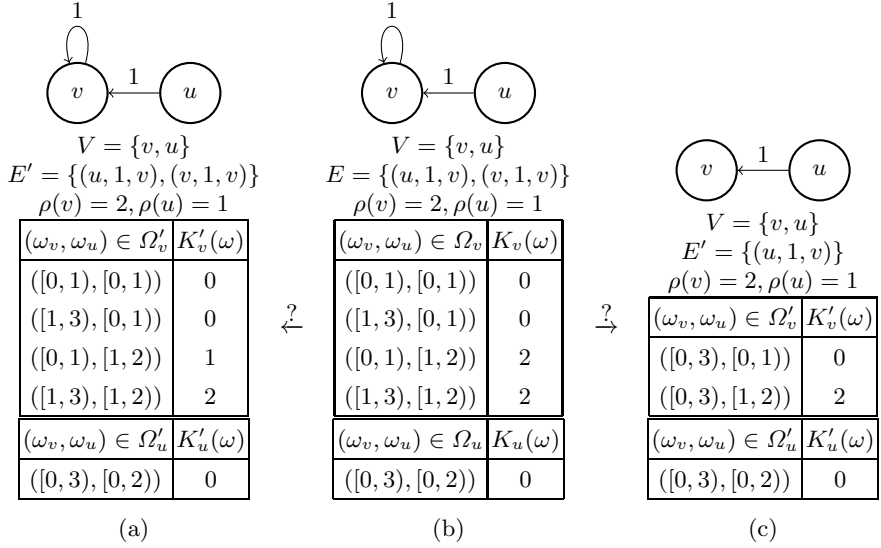
\begin{figure}[t]
	\centering
	\renewcommand{\arraystretch}{1.33}
	\begin{subfigure}[b]{0.3\textwidth}
		\centering
		\begin{tikzpicture}
	      \node[draw, circle, thick, minimum size=8mm] (1) at (0,0) {$v$};
	      \node[draw, circle, thick, minimum size=8mm] (2) at (1.5,0) {$u$};
	      \draw [->, in=105, out=75, loop] (1) to node[above] {$1$} (1) ;
	      \draw [->] (2) to node[above] {$1$} (1) ;
	    \end{tikzpicture}\\
    	$V = \{ v , u\} $\\
		$E' = \{ (u,1,v), (v,1,v) \}$\\
		$\rho(v) = 2, \rho(u)=1$ \\
    	\begin{tabular}{|c|c|}
    	\hline
    	$(\omega_v, \omega_u)  \in \Omega'_v$ & $K'_v(\omega)$  \\ \hline
    	$([0,1),[0,1))$ & 0 \\
    	$([1,3),[0,1))$ & 0 \\
      	$([0,1),[1,2))$ & 1 \\
    	$([1,3),[1,2))$ & 2 \\  	
    	\hline
    	\end{tabular}
    	\begin{tabular}{|c|c|}
    	\hline
    	$(\omega_v, \omega_u)  \in \Omega'_u$ & $K'_u(\omega)$  \\ \hline
    	$([0,3),[0,2))$ & 0\\
    	\hline
    	\end{tabular}
	    \caption{}
	    \label{fig:ambiguous_canonized}
    \end{subfigure}
	\begin{subfigure}[b]{0.025\textwidth}
	$\underleftarrow{?}$\vspace{80pt}
    \end{subfigure}
	\begin{subfigure}[b]{0.3\textwidth}
		\centering
		\begin{tikzpicture}
	      \node[draw, circle, thick, minimum size=8mm] (1) at (0,0) {$v$};
	      \node[draw, circle, thick, minimum size=8mm] (2) at (1.5,0) {$u$};
	      \draw [->, in=105, out=75, loop] (1) to node[above] {$1$} (1) ;
	      \draw [->] (2) to node[above] {$1$} (1) ;
	    \end{tikzpicture}\\
    	$V = \{ v , u\} $\\
		$E = \{ (u,1,v), (v,1,v) \}$\\
		$\rho(v) = 2, \rho(u)=1$ \\
    	\begin{tabular}{|c|c|}
    	\hline
    	$(\omega_v, \omega_u)  \in \Omega_v$ & $K_v(\omega)$  \\ \hline
    	$([0,1),[0,1))$ & 0 \\
    	$([1,3),[0,1))$ & 0 \\
      	$([0,1),[1,2))$ & 2 \\
    	$([1,3),[1,2))$ & 2 \\  	
    	\hline
    	\end{tabular}
    	\begin{tabular}{|c|c|}
    	\hline
    	$(\omega_v, \omega_u)  \in \Omega_u$ & $K_u(\omega)$  \\ \hline
    	$([0,3),[0,2))$ & 0\\
    	\hline
    	\end{tabular}
	    \caption{}
	    \label{fig:ambiguous_source}
    \end{subfigure}
	\begin{subfigure}[b]{0.025\textwidth}
	$\underrightarrow{?}$\vspace{80pt}
    \end{subfigure}    
    \begin{subfigure}[b]{0.3\textwidth}
		\centering
		\begin{tikzpicture}
	      \node[draw, circle, thick, minimum size=8mm] (1) at (0,0) {$v$};
	      \node[draw, circle, thick, minimum size=8mm] (2) at (1.5,0) {$u$};
	      \draw [->] (2) to node[above] {$1$} (1) ;
	    \end{tikzpicture}\\
    	$V = \{ v , u\} $\\
		$E' = \{ (u,1,v) \}$\\
		$\rho(v) = 2, \rho(u)=1$ \\
    	\begin{tabular}{|c|c|}
    	\hline
    	$(\omega_v, \omega_u) \in \Omega'_v$ & $K'_v(\omega)$ \\ \hline
   		$([0,3),[0,1))$ & 0 \\
    	$([0,3),[1,2))$ & 2 \\
    	\hline
    	\end{tabular}
    	\begin{tabular}{|c|c|}
    	\hline
    	$(\omega_v, \omega_u)  \in \Omega'_u$ & $K'_u(\omega)$  \\ \hline
    	$([0,3),[0,2))$ & 0 \\
    	\hline
    	\end{tabular}
	    \caption{}
	    \label{fig:ambiguous_reduced}
    \end{subfigure}
	\caption{Example of uncertain reduction. The network b) is in a non-canonical form. The canonization a) however makes the edge $(v,1,v)$ observable, even though it is superfluous as illustrated by c).}
	\label{fig:ambiguity}
\end{figure}

%% file: manuscript.bbl
\begin{thebibliography}{10}

\bibitem{bloniarz84a}
P.~A. Bloniarz, H.~B. Hunt, III, and D.~J. Rosenkrantz.
\newblock Algebraic structures with hard equivalence and minimization problems.
\newblock {\em Journal of the ACM}, 31(4):879--904, Sept. 1984.

\bibitem{dennunzio12a}
A.~Dennunzio, E.~Formenti, and J.~Provillard.
\newblock Non-uniform cellular automata: Classes, dynamics, and decidability.
\newblock {\em Information and Computation}, 215(0):32 -- 46, 2012.

\bibitem{gallet14a}
E.~Gallet, M.~Manceny, P.~Le~Gall, and P.~Ballarini.
\newblock An {LTL} model checking approach for biological parameter inference.
\newblock In {\em Formal Methods and Software Engineering}, volume 8829 of {\em
  Lecture Notes in Computer Science}, pages 155--170. Springer International
  Publishing, 2014.

\bibitem{gershenson04a}
C.~Gershenson.
\newblock Introduction to random boolean networks.
\newblock {\em arXiv preprint nlin/0408006}, 2004.

\bibitem{kauffman69a}
S.~Kauffman.
\newblock Metabolic stability and epigenesis in randomly constructed genetic
  nets.
\newblock {\em Journal of Theoretical Biology}, 22(3):437 -- 467, 1969.

\bibitem{streck12b}
H.~Klarner, A.~Streck, D.~\v{S}afr{\'a}nek, J.~Kol\v{c}\'ak, and H.~Siebert.
\newblock Parameter identification and model ranking of {Thomas} networks.
\newblock In {\em Computational Methods for Systems Biology}, volume 7605 of
  {\em Lecture Notes in Computer Science}, pages 207--226. Springer Berlin
  Heidelberg, 2012.

\bibitem{naldi11a}
A.~Naldi, E.~Remy, D.~Thieffry, and C.~Chaouiya.
\newblock Dynamically consistent reduction of logical regulatory graphs.
\newblock {\em Theoretical Computer Science}, 412(21):2207--2218, May 2011.

\bibitem{siebert08a}
H.~Siebert.
\newblock Local structure and behavior of {B}oolean bioregulatory networks.
\newblock In {\em Algebraic Biology}, volume 5147 of {\em Lecture Notes in
  Computer Science}, pages 185--199. Springer Berlin Heidelberg, 2008.

\bibitem{streck14a}
A.~Streck and H.~Siebert.
\newblock Equivalences in multi-valued asynchronous models of regulatory
  networks.
\newblock In {\em Cellular Automata}, volume 8751 of {\em Lecture Notes in
  Computer Science}, pages 571--575. Springer International Publishing, 2014.

\bibitem{thomas91a}
R.~Thomas.
\newblock Regulatory networks seen as asynchronous automata: A logical
  description.
\newblock {\em Journal of Theoretical Biology}, 153(1):1 -- 23, 1991.

\bibitem{thomas13a}
R.~Thomas.
\newblock Remarks on the respective roles of logical parameters and time delays
  in asynchronous logic: An homage to {El Houssine Snoussi}.
\newblock {\em Bulletin of Mathematical Biology}, 75(6):896--904, 2013.

\end{thebibliography}
